\documentclass[11pt, reqno]{amsart}
\usepackage{array}

\usepackage{float}
\usepackage{amsmath,latexsym}
\usepackage{amssymb}
\usepackage{geometry}
\usepackage{amsfonts}
\usepackage{hyperref}
\usepackage{color}

\usepackage{adjustbox}
\usepackage{booktabs}
\usepackage{tikz}
\usetikzlibrary{trees}

\usepackage{multirow}

\newcommand{\K}{\mathcal{K}}
\newcommand{\R}{\mathcal{R}}

\renewcommand{\Re}{\mathbb{R}}

\usepackage[numbers]{natbib}

\usepackage{setspace}
\allowdisplaybreaks

\let\origmaketitle\maketitle
\def\maketitle{
  \begingroup
  \def\uppercasenonmath##1{} 
  \let\MakeUppercase\relax 
  \origmaketitle
  \endgroup
}

\usepackage{graphicx}%
\usepackage{multirow}%
\usepackage{mathrsfs}%
\usepackage[title]{appendix}%
\usepackage{xcolor}%
\usepackage{textcomp}%
\usepackage{manyfoot}%
\usepackage{booktabs}%
\usepackage{listings}%
\usepackage{color}%
\usepackage{rotating}  
\usepackage{float}
\usepackage{adjustbox}

\usepackage{xcolor}

\title[]{\huge  Optimal probabilistic feature shifts for reclassification in tree ensembles }

\author[V. Blanco, A. Jap\'on, J. Puerto \MakeLowercase{and} P. Zhang]{{\large V\'ictor Blanco$^\dagger$, Alberto Jap\'on$^{\ddagger,\star}$, Justo Puerto$^\ddagger$, and Peter Zhang$^\star$}\medskip\\
$^\dagger$Institute of Mathematics (IMAG), Universidad de Granada\\
$^\ddagger$Institute of Mathematics (IMUS), Universidad de Sevilla\\
$^\star$ Heinz College of Information Systems and Public Policy, Carnegie Mellon University}
\date{}

\begin{document}

\maketitle
\vspace{-2em} 
\begin{center}
    \texttt{vblanco@ugr.es}, \texttt{ajaponsa@andrew.cmu.edu}, \texttt{puerto@us.es}, \texttt{pyzhang@cmu.edu}
\end{center}

\begin{abstract}
In this paper we provide a novel mathematical optimization based methodology to perturb the features of a given observation to be re-classified, by a tree ensemble classification rule, to a certain desired class. The method is based on these facts: the most viable changes for an observation to reach the desired class do not always coincide with the closest distance point (in the feature space) of the target class; individuals put effort on a few number of features to reach the desired class; and each individual is endowed with a probability to change each of its features to a given value,  which determines the overall probability of changing to the target class. Putting all together, we provide different methods to find the features where the individuals must exert effort to maximize the probability to reach the target class. Our method also allows us to rank the most important features in the tree-ensemble. The proposed methodology is tested on a real dataset, validating the proposal.

\end{abstract}

\section{Introduction}

Nowadays, machine learning models are used in a wide range of knowledge areas such as health care, transportation, finance, etc. Models are commonly required to fulfill certain characteristics that depend on the area of study, but generally there are two main features that stand out: precision and interpretability. For example, the classification models used to predict whether a credit card transaction is fraudulent or not, are as most accurate as possible, while it is less important to understand the rules that the model uses~\cite{london2019artificial}. On the other hand, there are other cases where besides accuracy, it is of vital importance to understand the classifier itself. For example when predicting a certain disease outcome, being able to interpret the disease progression implied by a model could be important to understand the nature of the disease, or how to mitigate it. Consequently, given the numerous advantages of using interpretable machine learning models for acquiring knowledge, they have been widely studied and have attracted much interest in recent years (see \cite{rudin2022interpretable}).\\ 

When thinking about interpretable machine learning, the first fundamental aspect to consider is the global interpretability of the model and how the variables of the problem interact. Linear and logistic regression are some of the preferred models because of their easy and manageable interpretation. In addition, different techniques, such as lasso or robust regression (see \cite{vidaurre2013survey,xu2008robust}) can be used to reduce the number of features to be interpreted within the models, which can significantly reduce the complexity of the task. On the other hand, decision tree-based models with univariate splits, or decision forests when they are of moderate size, are also very popular since the solutions can be reduced to a set of concatenated univariate decisions that can be easily followed and have a simple visual representation. In the latter, advances have been made in the field of mathematical optimization to derive optimal models that, smaller than their heuristic counterparts and reducing the number of univariate decisions to be interpreted, do not sacrifice accuracy when working with small to medium-sized data sets. In particular, in \cite{OCT} and \cite{OCF} the authors introduce Optimal Classification Trees (OCT) and Weighted Optimal Classification Forests, respectively, which have been proven to be useful in the design of interpretable classification tools. The multiple advantages of being able to interpret models globally are so widely spread that it can even be seen in \cite{craven1995extracting} and \cite{krishnan1999extracting} different ways of 
extracting decision trees from complex models such as neural networks, to try to understand how these models work internally.\\

On the other hand, there are cases in which, in addition to understanding the global functioning of the model, which may or may not be possible, it is of great interest to understand the prediction pathways for specific data points. For example, in binary classification, we will be interesting in understanding why a \emph{negative} label prediction is assigned to a particular observation. Specifically, we derive \emph{feasible} changes in the feature values of this observation to be predicted in the desired label. We will call the new perturbed observation the \emph{feature shifts} (FS) vector.  
 The analysis of the \emph{best} feature shifts for a given observation allows us to develop causal studies and open up predictive models to interventions results~\cite{gaggero2023sick,gaggero2024shutting}. The classical example in this type of studies is the case of an individual who applies for a mortgage loan in a bank and this is rejected (by certain classification rule). In this case, it would be beneficial for both the bank and the applicant to understand what the individual could change (assuming the classification tool stays the same) in the near future (e.g., if applicable, try to increase its salary), so that the loan would be granted, allowing the individual to buy the property and the bank to benefit from the loan. Nevertheless, to provide explanations on the reasons why an observation failed to be positively classified, it is clear, that the classification method must be explainable and interpretable. Furthermore, understanding how to \emph{repair} undesirable classification labels may have direct implications in the design of strategies, as for instance, in the development of public policies to palliate a disease, traffic accidents, suicides, etc. 

The study of FSs has been carried out in numerous models and applications, from the most interpretable to the most opaque or black-box models~\citep[see e.g.][]{guidotti2024counterfactual,verma2020counterfactual}. In our work, we focus on the study of FSs solutions for tree ensembles under mathematical optimization lens. Although this is not new, as it has already been studied in \cite{carrizosa2024generating} or \cite{parmentier2021optimal}, we take a different approach based on the feasibility of change for each individual, which can provide several advantages. Most of the proposals to derive FSs are distance-based models. They construct the most appropriate changes to perform on the features of one or more individuals to be classified as positive, using the closest-distance rule and the solution is the closest  point (with respect to a given distance) in the feature space to the observation/s. A drawback of this type of feature shifts is that the solutions may sometimes not be entirely satisfactory, since they do not necessarily reflect how feasible it is in practice for the observations to achieve the proposed solutions. For instance, if the feature shifts proposed in the solution for an observation is very different from the rest of the data, it may be perceived as impractical or impossible to achieve. On top of that, when using a distance-based approach, the solutions can sometimes fall very close to the decision boundary of the classifier, which makes them unstable to small perturbations of the classifier itself. Faced with these scenarios, in \cite{dutta2022robust} the authors propose robust alternatives to finding FSs. In addition, other methods have been proposed to obtain sparse solutions \cite{keane2020good,van2021interpretable}, maximize the connectivity from one point to another within a network, or look for close (in a general sense) solutions to the sample data \cite{dandl2020multi,forster2021capturing,kanamori2021distribution,mahajan2019preserving}. Although there are different metrics for comparing FSs, all seem to have advantages and disadvantages, and for the moment, a perfect and universal valid construction of FSs that fits all scenarios has not been found.

In this paper, we take a novel step forward in the construction of feasible and robust FSs by considering given probabilities of change in the values of  the features for an individual in the sample, i.e., the feasibility of changing the feature values. We present a mathematical optimization based approach that seeks to find the factors that maximize the probability that an  individual has to achieve a target class in the model. By downplaying the importance of the exact location of the feature shifted vector, we avoid the problems derived from choosing one metric or another to evaluate the solution.  On the contrary, our method puts all efforts on identifying which features the individual has to modify to achieve the objective with a higher degree of freedom. We do this by means of solving different Mixed Integer Non Linear Problems (MINLP). Moreover, it is important to note that, to reach a feature shift, the individual can work towards achieving the desired classification by focusing it effort into particular features (in the loan example, the individual might try to work overtime to get a salary increase, but not in reducing the loan amount, even if it was close to the initial feature value for this individual). The correct allocation of these efforts is also considered in our models as one of the main decision variables. We assume that the total \emph{effort} invested in changing its classification is upper bounded, which reflects the fact that, in practice, it is difficult (or even impossible) for the individual to improve the values of all the features. In addition, this allows our method to generate a new way to obtain smart feature importance rankings, that can provide advantages over other feature importance methods in tree ensembles. 

The rest of the paper is organized as follows. In section \ref{sec:2} we recall the general ideas of the FSs problem based on distances in tree ensembles (or forests) under the perspective of mathematical programming, and we set the notation used through the paper. In this section we also present our problem under study. Section \ref{sec:3} is devoted to present some valid MINLP formulations to solve our problem. In Section \ref{sec:4}, we present a case study, based on analyzing a problem for the prevention of obesity in some South American countries. Finally, some conclusions and further research on the topic are drawn in Section \ref{sec:5}.

\section{Preliminaries}\label{sec:2}

Throughout this paper, we consider that we are given an observation $x_0\in\mathbb{R}^d$ and a trained forest $\mathcal{F}$ comprising $R$ binary trees. We denote by $\R=\{1, \ldots, R\}$ the index set for these trees. 
This forest derives a decision rule $f:\mathbb{R}^d\rightarrow \K$, which assigns $x_0$ to a class $k_0 \in \K = \{1,\ldots, K\}$. In such a forest, each tree, $r\in \R$, predicts a class with weight $w^r\geq 0$. Gathering all the predictions amongst the trees, the decision rule follows a weighted majority voting strategy. Thus, if we consider $\Delta_k(x_0)$ to be the subset of trees that classify $x_0$ in class $k\in \K$, the previously mentioned class assignment occurs in case:
$$
k_0 = \arg\max_{k \in \K} \sum_{r\in\Delta_{k}(x_0)} w_r.
$$
Say the assignment to class $k_0$ is unsatisfactory, and that there exists a class $k^*\in \K$ which is the desired one for $x_0$. In this context, the solution to the FSs problem is to find a perturbation of $x_0$, $\delta(x_0)=x^* \in\mathbb{R}^d$, (by adequately altering the features of $x_0$), such that class $k^*$ is assigned to $x^*$ by $\mathcal{F}$, i.e.,
$$
k^* = \arg\max_{k \in \K} \sum_{ r\in\Delta_{k}(\delta(x_0))} w_r.
$$
In the literature, when this problem is approached  from the point of view of mathematical programming, 
it is often considered that there exists a cost function, $C(x_0,\delta(x_0))$, which represents the cost (usually a distance-based cost) of reaching $\delta(x_0)$ from $x_0$. Thus, the solution to the min-cost FSs problem is obtained by solving the problem
\begin{align}
\min_{x\in \Re^d} & \;\;  C(x_0,x)\nonumber\\
\mbox{ s.t. } & \sum_{r\in\Delta_{k^*}(x)} w^r \geq \sum_{r\in\Delta_k(x)} w^r, \forall k \in \K\backslash\{k^*\}.\label{cf:0}
\end{align}
The feature vector $\delta(x_0)$ solution of the problem above is known as the feature shifts (FS) of $x_0$.

\subsection{The distance-based approach}\label{subsect:2.1}

In distance-based approaches, $C(x_0,x)$ is induced by a certain  distance measure from $x_0$ to $x$, which is often easily implementable~\citep{blanco2014revisiting,blanco2024minimal}. Nevertheless, the same is not true for  constraints \eqref{cf:0} as it is required to properly model the assignment of the feature shifts to a suitable set of leaves in the forest, while defining the logistics and proper functioning of the forest itself. For a correct modeling, it must be taken into account that each tree in the ensemble is a set of branch nodes and leaf nodes, from which the observations follow a path from the root node to the leaves. Each branch node gives rise to two new nodes by means of a split until a leaf is reached, where a class is assigned to the observations that fall into it. The observations must verify the inequalities that define the splits in order to pass through them. 
In this paper, we consider the splits of the forest to be univariate, i.e., a single feature  intervenes to split each node.
Thus, these splits are of the form $x_{v(t)} \geq c_t$, so that at a node, $t$, the observation $x$ will follow the right branch if $x_{v(t)} \geq c_t$, and the left one if $x_{v(t)} < c_t$, where $v_{(t)}$ is the chosen variable, among the $d$ variables used to train the classifier appearing at node $t$, and $c_t$ is the threshold value. Thus, to formalize the modeling of the FSs problem, we will use the index sets and parameters described in Table \ref{tab:notation} for the data of the trained forest:
\begin{table}[h]
\begin{center}
\begin{tabular}{rp{10cm}}
\toprule
\multicolumn{2}{c}{Index Sets}\\
\midrule
$\tau^r$ & set of leaves of tree $r$.\\
$L_k^r$ & subset of leaves in tree $r$ whose output is class $k$.\\
$A_L^r(t)$& left ancestors of node $t$ in tree $r$ (nodes whose left branch has been followed on the path from the root node to $t$).\\
$A_R^r(t)$& right ancestors of node $t$ in tree $r$ (nodes whose right branch has been followed on the path from the root node to $t$).\\
\bottomrule
\toprule
\multicolumn{2}{c}{Parameters}\\
\midrule
$v_t^r$& feature involved in the split of node $t$ in tree $r$.\\
$c_t^r$& threshold value used in the split of node $t$ in tree $r$.\\
$w^r$& weight of tree $r$ in the forest.\\
\bottomrule
\end{tabular}
\end{center}
\caption{Parameters and index sets used in our models.\label{tab:notation}}
\end{table}

On the other hand, there are two sets of decision variables to consider, which comprise the location of the FSs solution in the space of variables, and the assignment of that solution to a set of leaves that gives the desired classification.  We define these variables as follows:
 \begin{itemize}
     \item $x \in \mathbb{R}^d$: feature shifts for $x_0$ in class $k^*$.
     \item $z_l^r \in \left\lbrace 0 , 1 \right\rbrace$: binary variable that takes value 1 if observation $x$ takes the path that leads to leaf $l$ in tree $r$.
\end{itemize}
Finally, for the sake of simplicity, and as it is usual in classification models, we assume that the features are normalized, taking values in $[0,1]$. Thus, one can assume that $c_t^r>\varepsilon$, $\forall t \in A_L^r(l) \cup A_R^r(l) , \  l \in\tau_l^r, \ r\in \R$, for a small enough $\varepsilon>0$.

With the above notation, the problem can be formulated as follows: 
\begin{align}
\displaystyle
\min  \;\; & C(x_0,x) & \nonumber\\
\mbox{s.t.} \nonumber \\
& x_{(v_t^r)} z_{l}^r \leq c_t^r - \varepsilon,& \forall t \in A_L^r(l),  l \in\tau_l^r, \ r\in \R, \label{c_1}\\
& (x_{(v_t^r)}-c_t^r) z_l^r  \geq 0,& \forall t \in A_R^r(l), l \in\tau_l^r, \ r\in\R, \label{c_2}\\
& \sum_{l\in\tau_l^r} z_l^r =1, & \forall r\in\R,\label{c_3}\\
& \sum_{r\in \R}\sum_{l\in L_{k^*}^r} w^r z_l^r \geq \sum_{r\in\R}\sum_{l\in L_{k}^r} w^r z_l^r, & \forall k \in \K\backslash\{k^*\},\label{c_4}\\
& z_{l}^r \in \{0,1\}, &\forall l \in \tau_l^r, r \in \R,\\
& x \in \mathbb{R}^d.&
\end{align}
Constraints \eqref{c_1} and \eqref{c_2} guarantee that if observation $x$ is in leaf $l$ of tree $r$ ($z_l^r = 1$), then observation $x$ must satisfy all the inequalities in the path from the root node to $l$ throughout all its ancestors. Furthermore, note that these constraints are redundant when $z_l^r = 0$. These constraints together with Constraints \eqref{c_3} force the feature perturbation, $x$, to be in one, and only one, leaf on every single tree in the forest. Finally, Constraints \eqref{c_4} assure that the majority voting strategy is carried out properly and thus the prediction for $x$ is the desired class $k^*$.

Note that the model above is a Mixed Integer Non Linear optimization problem both because of the the objective function and constraints \eqref{c_1} and \eqref{c_2}. However, as already mentioned, if $C$ is induced by a $\ell_p$ or a polyhedral norm, it can be efficiently represented by a set of second order cone constraints~\cite{blanco2014revisiting,blanco2024minimal}, and the constraints can be linearized by using big-M constraints, resulting in a Mixed Integer Second Order Cone Optimization (MISOCO) problem.

\subsection{Our approach}

We propose a different approach than considering a cost function based on distances to be minimized. We set that when observation $x_0$ obtains an undesired classification by the forest, it will take some time for $x_0$ to go through the classifier again, and in this time $x_0$ has some probabilities of changing the values of its variables, and besides, $x_0$ can also try to intervene so that these changes are favorable for the desired classification. At the same time, we assume that the ability of $x_0$ to intervene favorably in this task is limited. As usual in real-life applications, it is impractical to assume that all the features in $x_0$ can be altered to reach the target class $k^*$. Thus, $x_0$ is provided with a threshold value, representing the maximum total allowed effort towards the target class. Our problem is then focused on  finding the feature perturbations for $x_0$ by optimally allocating these limited efforts, in the sense of what is most likely to occur in the future.

The motivation behind our problem is two-fold. First, we must consider that the ``closest'' feature shifted solution is not necessarily the most likely or feasible one when the individual is re-evaluated by the forest. For instance, in datasets with categorical features (encoded by means of binary variables), norm-based distances between the original $x_0$ and a feasible feature shifts  $x$ may be large, whereas, in practice, changing a single binary feature could be highly probable. Additionally, many other problems can arise when considering distances with variables of magnitudes that are not easily comparable. We avoid this situation by considering a probabilistic approach. Secondly, when a forest is used to obtain an importance ranking of variables, in general, these are obtained through the splits that intervene in the trees based on their classification capacity for the sample data with which the classifier has been trained. However, there are features that an individual is not allowed to change (such as age, gender, etc) or that its probability to change is small (occupation, geographic district of home address, etc) which, having been used in the classifier, have an impact on the ranking of variables. Therefore, if these variable rankings are used to design strategies for changing an individual's classification into the desired class, the strategies may be suboptimal, as they can be misled by the effects of these fixed variables in the ranking. On the contrary, by using our method, we can construct feature importance rankings based on observing where the efforts of our feature shifts  solutions are \emph{optimally} allocated among the features, thus avoiding to involve the effect of fixed variables (having 0 probability of change), and resulting in the design of more efficient and realistic strategies.

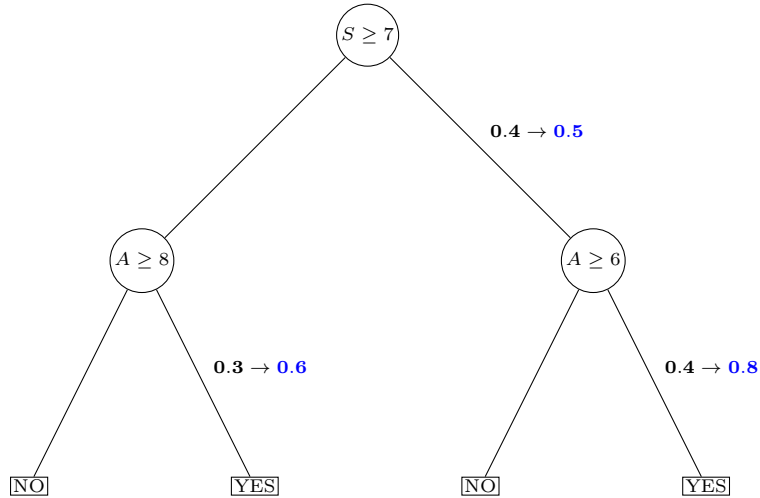
\begin{figure}[h]
    \centering

\begin{tikzpicture}[scale=2, level distance=1.5cm,
level 1/.style={sibling distance=3cm},
level 2/.style={sibling distance=1.5cm}]

\node[circle,draw, inner sep=1pt] (Root) {\tiny $S\geq 7$}
    child {
    node[circle,draw, inner sep=1pt](left1) {\tiny $A \geq 8$} 
    child { node[rectangle, draw, inner sep=1pt] {\tiny NO} }
    child { node[rectangle, draw, inner sep=1pt](right2) {\tiny YES} } 
}
child {
    node[circle,draw, inner sep=1pt](right1) {\tiny $A\geq 6$}
    child { node[rectangle, draw, inner sep=1pt] {\tiny NO} }
    child  { node[rectangle, draw, inner sep=1pt](right3) {\tiny YES} }
};
\path (Root) -- node[above right] {\tiny $\mathbf{0.4 \rightarrow {\color{blue} 0.5}}$} (right1);
\path (left1) -- node[above right] {\tiny $\mathbf{0.3 \rightarrow {\color{blue} 0.6}}$} (right2);
\path (right1) -- node[above right] {\tiny $\mathbf{0.4 \rightarrow {\color{blue} 0.8}}$} (right3);

\end{tikzpicture}
    \caption{Firefighter example}
    \label{fig:1}
\end{figure}

We end this section illustrating our approach in a forest of size 1, i.e. a tree. Suppose that an individual $x_0$ desires to become a firefighter. To do so, the individual must pass an entrance exam at an academy, which did not happen on the first attempt. The test consists of two scores resulting from a strength test (S) and an aerobic test (A), which are combined to determine whether an individual is accepted using the tree diagram shown in Figure \ref{fig:1}. The individual will have to score higher than 8 in test A when they have not achieved a minimum of 7 in test S, or will have to score higher than 6 in test A otherwise to be admitted. We assume that the individual has failed the test on the first try, e.g., having already been classified into a ``No'' leaf via one specific path. For the second try that is happening in the future, the outcome is random (each path from the root to a leaf has a certain probability of happening). In Figure \ref{fig:1} we can see in boldfaced black the probabilities of the individual to follow the right branch (in the next attempt to happen in the future) of the nodes when no effort is applied to the node variable. In color blue we indicate the probabilities when effort is applied. We always assume either that these probabilities are known or that they can be estimated, and they can be different for each individual. In this example, we further set that the individual has limited time and effort: they can only invest effort into improving one of the two tests. Being a toy example, we can easily check that if the individual puts effort on S, we would have a ``Yes'' (the left one) with probability $(1-0.5)\times0.3 = 0.15$, and another ``Yes'' (the right one) with probability $0.5\times0.4=0.2$. If the effort is put on A, we would have a ``Yes'' with probability $(1-0.4)\times0.6=0.36$, and the other ``Yes'' with probability $0.4\times0.8=0.32$ (observe that the combinatorial complexity of the problem increases significantly with the number of trees and when the majority vote decision rule is applied). At this point, there are now different strategies for deciding where to invest the effort. For instance, one option could be to look for the best-case scenario and to obtain the feature shifts $x$ associated with the highest probability path. Following this approach, the solution would be to put effort into test A, obtaining a best-path probability of 0.36, where the solution $x$ would be defined by obtaining scores such that $S < 7$ and $A\geq 8$. Alternatively, one can also consider a worst-case approach -- the individual picks the effort investment that results in the highest worst-case path probability. In this case the individual should also invest into test A, since it leads to the best worst-case probability of 0.32.
 In the next section, we will introduce different optimization models that  are capable of producing a solution under any user-specified risk level, ranging from the most conservative (maximizing worst-case path probability towards ``Yes''), to risk neutral (maximizing median-case path probability), to the most risk-seeking (maximizing the best-case path probability).

\section{A Mathematical optimization model to construct optimal feature shifts maximizing its occurrence probability}\label{sec:3}

We begin this section by presenting the solution to the problem of finding the feature shifts that is most likely to occur when the individual is again re-evaluated. We will assume given probabilities of taking the right (and left) branch in the future, when extra effort is applied and when no extra effort is applied, for every node of every tree in the forest. After this, we will also present the solution of different, more robust versions of the problem. For this purpose, in addition to the notation and elements discussed in the previous section, we need to introduce new ones and delve into some concepts. First, we have defined the total effort of an individual as a bounded resource (such as money or time), therefore we will denote by $\eta>0$ the maximum allowed amount of effort for the individual. 

We consider that this effort is measurable in units of effort, and that it can be distributed unequally among the set of variables, with $E$ being the maximum number of units of effort that an individual can put into a given feature. 

We denote by $p_{lt}^r$ the probability of taking in node $t$ the branch that leads to leaf $l \in \tau_l^r$ of tree $r$. Each unit of effort that is invested in improving the feature involved in such a node has a direct impact on that probability. We denote by $p_{lte}^{r}$ the probability  of taking in node $t$ the branch that leads to leaf $l$  of tree $r$ when $e \leq E$ units of effort are used.

Although our model can be extended to a general case for a classification problem with $K$ classes and having the trees of the forest weighted by weights, $w_r$, for the sake of simplicity, from now on we will consider a two-class (enconded as $0$ and $1$) instance and the trees weighted equally in the ensemble. 

\subsection{Most likely path problem}

To define the objective of the problem of finding the most probable feature shifted solution in the forest, we must first note that a binary classification forest with $R$ equally weighted trees gives a classification in class $k^*$ for $x\in \mathbb{R}^d$ if at least $R/2$ of the trees classify $x$ in $k^*$. In order to know which vote each tree gives, it will then be important to model on which leaf of the tree our candidate for solution $x$ is located. This decision is provided by the $z$-variables already introduced in Subsection \ref{subsect:2.1}. On the other hand, since the trees are independent, the probability of a given vote occurring in the forest can be decomposed into the product of the probabilities of the trees that give that vote. Once a distribution of effort is made, we would have defined the probability that observation $x$ reached leaf $l$ on tree $r$. Since we only need $R/2$ of the trees to give a favorable vote to classify in $k^*$, we now define as \emph{essential} trees those $R/2$ that give a vote in $k^*$ with the highest probability. We define the variables $\theta_l^r\in \left[0,1\right]$ as the probability of reaching leaf $l$ in tree $r\in \R$, given an effort distribution, when $r$ is an essential tree. It takes value 1 when the tree is not essential. The probability to reach the feature shifts $x$ from $x_0$ can be written then as:
$$
\prod_{r\in \R} \sum_{l\in\tau_l^r}z_{l}^{r} \theta_l^r, 
$$
This is the proposed objective function to be maximized by  our model. By monotonicity of the $\log$ function, maximizing this probability will be equivalent to maximize its logarithm. Thus, our objective function for this model will be:
\begin{align}
\displaystyle
\max_{z, \theta}  \;\; & \log \left( \prod_{r\in \R} \sum_{l\in\tau_l^r}z_{l}^{r} \theta_l^r   \right) = \sum_{r\in \R} \log\left(\sum_{l\in\tau_l^r}z_{l}^{r} \theta_l^r   \right).\nonumber  
\end{align}
According to this, as $x$ can only be located in one leaf of each tree, we would look for a feasible valid allocation in the $R/2$ essential trees that will maximize the probability of the desired classification.  

Having established the objective function, we then need to add constraints to the problem to ensure the correct definition of the variables, the functioning of the forest and to guarantee that we obtain a valid feature shifted solution. To begin with, it should be noted that the constraints \eqref{c_1}-\eqref{c_4} discussed in Section \ref{sec:2} for the distance-based feature shifts problem must also be added in our problem. This will guarantee that the observation $x$ is in only one of the leaves of each tree, and that the splits of the inequalities of the nodes leading to that leaf are verified. Furthermore, it will also guarantee that $x^*$ is a valid feature shifts solution as there will be more trees giving a vote in $k^*$ than in the opposite class.

The next element to consider in our model will be the correct definition of the effort and its related variables. To do this, we incorporate the following set of variables to our model:
$$
\beta_{te}^r = \begin{cases}
    1 & \mbox{if $e$ units of effort are applied to improve the value of the}\\
    & \mbox{feature involved in the split of node $t$ in tree $r$},\\
    0 & \mbox{otherwise.}
\end{cases}
$$
for all $e =0,\ldots,E, \ t \in\tau^r,  r \in \R$.

With this definition, we ensure that in all variables we apply an effort from $0$ units--no effort, to a maximum of $E$ units--maximum allowed effort). Furthermore, using the $v$-parameters already described above, we must assure that if a feature appears in more than one split of the whole forest, the associated variables $\beta_{te}^r$ must have the same value in all occurrences. These two requirements are assured by means of the following constraints:
\begin{align}
& \sum_{e=0}^E \beta_{te}^r = 1,  & \forall  t\in\tau^r ,\ r\in\R, \label{c_5} \\
& \beta_{te}^r = \beta_{t'e}^{r'}, \text{\ if} \ v_t^r = v_{t'}^{r'} ,  &  \  \forall e =0,\ldots,E, \ t,\ t' \in\tau^r, \ r, \ r'\in \R.  \label{c_6}
\end{align}
In addition, one must also ensure that no more than $\eta$ units of total effort are applied to construct the feature shifted vector. If we denote by $b_t^r$ the total number of times that variable involved in node $t$ of tree $r$ appears in the forest, this can be  done by adding the following constraint to the problem:
\begin{align}
& \sum_{r\in \R}\sum_{t\in\tau^r} \sum_{e=0}^E\frac{e\beta_{te}^r}{b_t^r} \leq \eta. \label{c_7} &
\end{align}
Finally, to adequately define the $\theta_l^r$ variables (probability of reaching leaf $l$ in tree $r$), we must note that these variables are maximized in the objective function, and that the task is therefore to find the correct upper bounds for them. As mentioned above, when these variables are associated with non-essential trees, they must be allowed to take value 1. A tree can be non-essential, either because the vote it gives is in the undesired class, or because despite giving a favorable vote, there are $R/2+1$ trees that give a favorable vote with a higher probability. We will then add some extra binary variables and indicators to the problem to control these two scenarios:
$$
\gamma^r = \begin{cases}
    1 & \mbox{if tree $r$ gives a positive vote that is not essential to}\\
    & \mbox{obtain the desired classification,}\\
    0 & \mbox{otherwise}
\end{cases} \quad \forall r \in \R.
$$
To properly define these variables we need to consider two factors. First, we have to force the variables to take value $0$ when the trees do not assign the desired class to $x$. If we denote by $\iota_l^t$ the predicted class at leaf $l$ of tree $r$, this condition can be imposed with the following constraints:
\begin{align}
& \gamma^r \leq \sum_{l\in\tau_l^r : \iota_l^r=k^*} z_l^r,  & \forall r\in \R. \label{c_8}
\end{align}
Additionally, we must ensure that we allow these variables to take value $1$ only in the appropriate trees, and for this we have to regulate the total number of variables that take value $1$. Since $(R/2+1)$ of the trees are essential, we will have to, at all times, control how many trees give a favorable vote to $x$, and subtract from this number $(R/2+1)$, as this would be the number of non-essential trees that vote for the desired class. We do this by adding to the problem the following constraint:
\begin{align}
&\sum_{r\in\R} \gamma^r \leq \sum_{r\in\R} \sum_{l\in\tau_l^r : \iota_l^r=k^*} z_l^r-\left(\frac{R}{2}+1\right). & \label{c_9}
\end{align}

Furthermore, when modeling non-essential trees we need to avoid those leaves of trees that do not predict the target class for the feature shifts. We consider the following binary values:
$$
\phi_l^r = \begin{cases}
    1 & \mbox{if the predicted class in leaf $l$ of tree $r$}\\
     & \mbox{is not the desired class},\\
    0 & \mbox{otherwise}, 
\end{cases}
\quad \forall l \in \tau^r, \forall r \in \R.
$$
Note that $\phi_l^r$ can be easily calculated in preprocess (this is not a variable of our model, but a parameter).

With all the above, we can assure the correct definition of the $\theta_l^r$ variables with the following constraints:
\begin{align}
& \theta_l^r  \leq \prod_{t\in A(l)} \left( p_{lt}^r + \sum_{e=0}^{ E}\beta_{te}^r(p_{lte}^{r} - p_{lt}^r) \right) + \gamma^r + \phi_l^r, & \forall l \in\tau_l^r,r\in \R. \label{c_10}    
\end{align}
In these constraints, the first term on the right-hand-side refers to the probability of reaching leaf $l$ in tree $r$, which is equal to the product at the ancestor nodes of $l$, $A(l)$, of the expression $p_{lt}^r + \sum_{e=0}^{ E}\beta_{te}^r(p_{lte}^{r} - p_{lt}^r)$. In this way, if no effort is applied to any of the variables of the nodes arriving at leaf $l$ ($\beta_{te}^r = 0$ for all the nodes in the branch), we would simply be calculating the product of the $p_{lt}^r$, and when e units of effort were applied to any node ($\beta_{te}^r = 1$), the probability of these nodes would be updated in the product to $p_{lte}^r$. In addition, if a leaf does not give the desired classification, the first term referring to the probabilities would have no effect to $\theta_l^r$ as $\phi_l^r$ would take value $1$. A similar case would occur if the positive vote is given by a non-essential tree, since $\gamma^r$ would take value $1$ as well. In fact, it is these constraints that determine which trees give the essential votes and which give the non-essential votes. Note that by constraints \eqref{c_4}, at least $(R/2+1)$ trees are forced to give a positive vote for the feature shifts solution to be valid, but in case of having more, by maximizing $\theta_l^r$ in the objective function, the solution will implicitly convert into non-essential ($\gamma^r =1$) those that do so with less probability.\\

Gathering all together, the problem is formulated as follows:
\begin{align*}
\displaystyle
\max   & \;\;\sum_{r\in \R} \log \left(\sum_{l\in\tau_l^r}z_{l}^{r} \theta_l^r   \right) &  \tag{${\rm Max-path \ }$} \\ 
\mbox{s.t. } & \eqref{c_1}-\eqref{c_10},\\
& \theta_l^r \in [0,1], & \forall l \in \tau^r, r \in \R,\\
& z_l^r \in \{0,1\}, & \forall l \in \tau^r, r \in \R,\\
& \gamma^r\in \{0,1\}, &\forall r\in \R,\\
& \beta_{te}^r \in \{0,1\}, & \forall e=0, \ldots, E, t\in \tau^r, r\in \R,\\
& x\in \Re^d.
\end{align*}

This is a MINLP problem due to the logarithm function appearing in the objective function, and the product of variables. However, the products are easily linearizable with the help of some auxiliary variables, and the logarithm function can be accurately (outer) approximated by piecewise linear functions.

\subsection{Robust versions of the problem}

In the previous section we have presented the problem of distributing a certain amount of effort so as to obtain the most likely feature shifts solution to occur in the future. This problem can yield very good results in certain applications, but like any methodology, it has some advantages and disadvantages. In particular, this method can work well when there is an allocation of effort that gives a high probability to our feature shifts solution. However, as the number of trees increases, as well as their depth, the number of leaf combinations on the trees reporting the desired classification through the voting strategy increases exponentially, and therefore, generally, the probability of occurrence of exactly one of these combinations tends to be small. This fact may call into question the sense of studying the distribution of efforts in order to maximize a single scenario, and this naturally opens the door to study more robust versions of the problem.\\

An alternative option would be to find the $R/2+1$ essential trees, but instead of maximizing their path to the leaf with the highest probability of the desired class, maximize the worst-case scenario, i.e., the one with the smallest  probability. This approach will distribute the effort to maximize the probability of the valid feature shifts solution least likely to occur, and by maximizing this, we would be guaranteed to be distributing effort with some impact on all possible voting combinations. In this problem, instead of maximizing the probability of the best scenario, we maximize the minimum probability, i.e., the function:
$$
\prod_{r\in \R} \left(\min_{l \in \tau^r} \theta_l^r\right)
$$
Thus, analogously to the previous model, the problem can then be formulated as follows:
\begin{align*}
\max & \;\;\sum _{r\in \R} \log\left(\min_{l \in \tau^r} \theta_l^r\right)\\
\mbox{s.t. } 
&\eqref{c_1}-\eqref{c_10},\\
& \theta_l^r \in [0,1], & \forall l \in \tau^r, r \in \R,\\
& z_l^r \in \{0,1\}, & \forall l \in \tau^r, r \in \R,\\
& \gamma^r\in \{0,1\}, &\forall r\in \R,\\
& \beta_{te}^r \in \{0,1\}, & \forall e=0, \ldots, E, t\in \tau^r, r\in \R,\\
& x\in \Re^d.
\end{align*}

As usual in this type of max-min problems, the minimum can be linearized by  introducing a new continuous auxiliary variable $y^r \in [0,1]$, that will take the value of the minimum of the $\theta_l^r$ variables amongst the leaves of the tree, resulting in the model:
\begin{align*}
\max & \;\;\sum _{r\in \R} \log y_r \tag{${\rm Min-path \ }$}\\
\mbox{s.t. } &\eqref{c_1}-\eqref{c_10}\\
& y^r \leq \theta_l^r , & \forall l \in \tau_l^r,   r\in \R,\\
& y^r \geq \theta_l^r - (1-z_l^r), & \forall l \in \tau_l^r, r\in \R,\\
& \theta_l^r \in [0,1], & \forall l \in \tau^r, r \in \R,\\
& z_l^r \in \{0,1\}, & \forall l \in \tau^r, r \in \R,\\
& \gamma^r\in \{0,1\}, &\forall r\in \R,\\
& \beta_{te}^r \in \{0,1\}, & \forall e=0, \ldots, E, t\in \tau^r, r\in \R,\\
& x\in \Re^d.
\end{align*}

Note that in case $r$ is a non-essential tree $y^r$ takes value $1$.

Despite gaining robustness against the distribution of efforts to maximize a single scenario, this method might  have some inconvenience. Although maximizing the probability of the least likely solution implicitly has an impact on all other solutions, in the objective function we are still optimizing only one path, and it tends to have a very small probability as the size of the forest increases. In this way, since a logarithm is involved in the objective function, numerical instabilities can occur when solving the problem and even report infeasibility due to negative values when the objective value approaches zero. In addition, one could argue that despite obtaining a robust solution, one might be focusing too much on optimizing for combinations that are very unlikely to occur, which despite giving a favorable classification, require paths in the trees that are impractical for the feature shifts solution. This is why we develop a family of models in between the two models already presented.

The cumulative probability of the least likely paths occurring in the trees can be very small. Thus, by setting a minimum likelihood threshold, we propose a model where instead of optimizing the least likely path, we optimize the least likely path above this threshold, thus decreasing the effect of artifacts that may come from very unlikely paths, while obtaining a robust solution.

We apply an approach based on the concept of $k$-sum optimization (also named in portfolio optimization \textit{Conditional-Value-at-Risk (CVaR)}). CVaR has been proven to be useful in different applications for deriving robust solutions for problems under risk or uncertainty \cite{fernandez2019new,filippi2020conditional,li2023risk,li2024risk,puerto2016}.

For a given probability threshold $\mu$, the goal is to find the feature shifts for $x_0$ such that the joint probability of the $\kappa$ less probable paths is maximized for a given $\kappa$ value. It would allow not only to consider the worst-case situation but a wide family of situations.

To incorporate that into a mathematical optimization model, for each tree $r\in R$ we need to sort the probabilities $\theta_{l}^r$, for all $l\in \tau^r$ (say $1, 2, 3, \ldots, |\tau^r|$) in non decreasing order, i.e., to find a permutation of the leaves, $\sigma_r$ that verifies:
\begin{equation}
\theta_{\sigma_r(1)}^r \leq \theta_{\sigma_r(2)}^r \leq \cdots \leq \theta_{\sigma_r(|\tau^r|)}^r, \quad \forall r\in \R.\label{sorting}
\end{equation}
This permutation will be modeled with the following assignment variables:
$$
\alpha_{il}^r = \begin{cases}
    1 & \mbox{if $\sigma_r(l)=i$}\\
    0 & \mbox{otherwise}
\end{cases}, \quad \forall i=1, \ldots, |\tau_r|,  r\in \R.
$$
The adequate definition of these variables is assured with the following constraints.

First, each leave must be sorted in a unique position in the sorting sequence of probabilites, and each position must be assigned to exactly one leaf of the tree:
\begin{align}
    &\sum_{i\in\tau_l^r} \alpha_{il}^r =1, & \forall l \in\tau_l,\  r\in\R, \label{assign:1}\\
&\sum_{l\in\tau_l^r} \alpha_{il}^r =1,& \forall i \in\tau_l,\  r\in \R.\label{assign:2}
\end{align}

Second, the sorting must verify \eqref{sorting}:
\begin{align}
 \sum_{l\in \tau^r} \alpha_{il}^r  \theta_l^r \leq \sum_{l\in \tau^r}  \alpha_{(i+1)l}^r  \theta_l^r, & \; \forall i=1, \ldots, |\tau^r|-1,  r\in \R.\label{sorting:2} 
\end{align}
that is, the leaf sorted in position $i$ must have probability smaller than the one in position $i+1$.

Note that the above constraint is nonlinear, although it can be linearized using standard McCormick envelopes.

Finally, we impose that the cumulative probability of the $\kappa-1$ smallest leaves of tree $r$ is at least $\mu$:
\begin{align}
\sum_{i=1}^{\kappa-1} \sum_{l \in \tau^r} \alpha_{il}^r \theta_{l}^r \geq \mu, & \;\; \forall r\in \R.\label{threshold}
\end{align}

The problem can then be formulated as follows:
\begin{align*}
\max & \;\;\sum_{r\in \R} \log\left(\sum_{l \in \tau^r} \alpha_{\kappa l}^r \theta_{l}^r\right) \tag{${\rm \kappa-path \ }$}\\
\mbox{s.t. } &\eqref{c_1}-\eqref{c_10},\\
& \eqref{assign:1}-\eqref{threshold},\\
& \theta_l^r \in [0,1], & \forall l \in \tau^r, r \in \R,\\
& z_l^r \in \{0,1\}, & \forall l \in \tau^r, r \in \R,\\
& \gamma^r\in \{0,1\}, &\forall r\in \R,\\
& \beta_{te}^r \in \{0,1\}, & \forall e=0, \ldots, E, t\in \tau^r, r\in \R,\\
& \alpha_{il}^r \in \{0,1\}, & \forall i=1, \ldots, |\tau^r|, l\in \tau_l^r, r\in \R,\\
& x\in \Re^d.
\end{align*}

This formulation after linearizing products results in the following mixed integer nonlinear problem.

\begin{align*}
\max & \;\;\log\prod_{r\in \R} u_{\kappa}^r \tag{${\rm \kappa-path \ }$}\\
\mbox{s.t. } 
& u_i^r \leq \theta_l^r + (1-\alpha_{il}^r), &\forall i, l \in\tau_l,\  r\in \R, \\
& \sum_{i\in\tau_l^r} \alpha_{il}^r =1, &\forall l \in\tau_l,\  r\in \R, \\
& \sum_{l\in\tau_l^r} \alpha_{il}^r =1,&\forall i \in\tau_l,\  r\in \R,\\
& u_i^r \le u_{i+1}^r,  & \forall i\in\tau_l,\forall r\in \R, \\
& \sum_{i=1}^{\kappa-1} u_i^r \ge \mu,&\forall r\in\R,\\
 &\eqref{c_1}-\eqref{c_10},\\
& \theta_l^r \in [0,1], & \forall l \in \tau^r, r \in \R,\\
 & z_l^r \in \{0,1\}, & \forall l \in \tau^r, r \in \R,\\
& \gamma^r\in \{0,1\}, &\forall r\in \R,\\
& \beta_{te}^r \in \{0,1\}, & \forall e=1, \ldots, E, t\in \tau^r, r\in \R,\\
& \alpha_{il}^r \in \{0,1\}, &\forall i=1, \ldots, |\tau^r|, l\in \tau_l^r, r\in \R,\\
& x\in \Re^d.
\end{align*}

\section{Case study: probabilistic feature shifts in obesity classification.}\label{sec:4}

In this section, we analyze a case study in which we demonstrate an important practical application derived from the use of our models. Specifically, we focus on the study of a classification problem to detect obesity, and we conduct our research based on the real data collected for the work developed in \cite{de2019obesity}, which can be found freely available on the Kaggle website. The aim of the work is to obtain different importance rankings of variables from a random forest classifier, and to evaluate the favorable/unfavorable effect that modifying the values of the most important variables amongst the ranking would have on the prediction of obesity. We create these rankings by properly using the solutions of our models, and compare them with those provided by the \texttt{feature\_importances\_} method of the \texttt{scikit-learn} python library as well as by a set of random rankings. Once the rankings are created, we run simulations to evaluate the robustness of the rankings. For this, we randomly perturb the initial data to create possible future scenarios, and evaluate these data with the classifier to see if there have been any changes with respect to the original obesity classification.\\

According to this, the problem to be addressed is to find feature shifts solutions for individuals that have been classified as obese, in order to detect the best effort allocation to modify their classification. The hypothesis of the experiment is that, if we put effort into modifying favorably the most important features with respect to the target of reclassifying individuals as healthy, we should obtain a greater number of people that can be reclassified as healthy than if we put no effort. For simplicity, the entire experiment is conducted with a single level of effort, therefore the decisions are binary, either to put or not to put effort into the available variables.\\

This section is divided in three parts. We start by introducing the data and the random forest classifier. After this, in order to be able to use our models, we describe and apply a methodology for estimating the probability of change at the nodes of the trees. We conclude by reporting the different rankings and results obtained in our  simulations.

\subsection{Data and forest}

The dataset we use in this case study was collected from  \cite{de2019obesity}. The dataset includes a total of $2111$ individuals between $14$ and $61$ years old from the countries of Mexico, Peru and Colombia, and were collected through an anonymous online survey. The dataset has  $16$ predictors and one target variable that tracks a person's weight status. This target variable has $7$ categories: underweight, normal weight, two levels of overweight, and three levels of obesity. Given that the interest of our study is focused on analyzing the factors that have the greatest impact on obesity, we have eliminated from the data the individuals with insufficient weight, and we have re-coded the target variable as a binary one, grouping the three types of obesity into one class, and grouping the normal weight and the two levels of overweight in another. On the other hand, we also eliminated two of the predictor variables, height and weight, as their combination was strongly related to the target variable but the information provided was not of our interest for the experiment. The resulting variables used are shown in Table \ref{tab:1}, where their abbreviation and name, description and univariate impact with respect to the target variable are listed. For continuous variables, this impact is positive, $+$, (resp negative)  if a higher increase in the variable is associated with a lower (resp higher) obesity rate. For binary variables, the impact is positive if belonging to class $1$ (resp. $0$) is associated with a lower (resp. higher) obesity rate. Finally, we observed that some of the classes of the predictor variables CALC and PT were extremely unbalanced. To address this, we removed from the sample the only individual with the highest daily alcohol consumption, and re-coded the original transport variable by a new binary one encoding the use of private or public transport.\\

\begin{table}[]
    \centering
    \begin{adjustbox}{max width=\textwidth}
    \begin{tabular}{l c c}
    \toprule
        {\bf Name} & {\bf Type} & {\bf Change Impact on Obesity} \\ \midrule
        Gender & Binary  &   \\
        Age & Continuous &  \\
        Family history with overweight (FHWO) & Binary &  \\
        Frequent consumption of high caloric food (FAVC) & Binary & Class 0 preferred \\
        Frequency of consumption of vegetables (FCVC) & Continuous  & - \\
        Number of main meals (NCP) & Continuous & + \\
        Consumption of food between meals (CAEC) & Continuous & + \\ 
        Smoke & Binary & Class 0 preferred \\
        Consumption of water daily (CH2O) & Continuous & - \\ 
        Calories consumption monitoring (SCC) & Binary & Class 1 preferred \\
        Physical activity frequency (FAF) & Continuous & - \\
        Time using technology devices (TUE)  & Continuous & + \\ 
        Consumption of alcohol (CALC) & Continuous & + \\
        Private Transport (PT) & Binary & Class 0 preferred \\
        \bottomrule
    \end{tabular}
    \end{adjustbox}
    \caption{Features in the dataset.}
    \label{tab:1}
\end{table}

After having preprocessed the data, we randomly partitioned the dataset into a training set (consisting on $2/3$ of the data), from which we trained the random forest and applied our models in order to obtain the variable rankings, and a test set (consisting on $1/3$ of the data) where we ran the simulations to evaluate the different rankings. Given the large number of individuals on which we were going to train our mathematical programming models, we needed the random forest on which the experiment is based to be of a manageable size for us. After several preliminary trials, we observed that our models provided accurate solutions with a time-limit of 300 seconds when the trees in the forest had a depth less than $7$. With these parameters, we used the \texttt{scikit-learn} python library to train the random forest. We concluded to use a random forest of 25 trees of depth equal to $5$, obtaining an accuracy rate of $85\%$, with an accuracy rate of $93\%$ in the obese group of the sample. Increasing the depth of the trees to $6$, or the number of trees to $100$ had a marginal impact on the classification rate. We therefore decided to use the smallest model to reduce the computational cost.\\

\begin{figure}[h]
    \centering
    \includegraphics[width=0.75\linewidth]{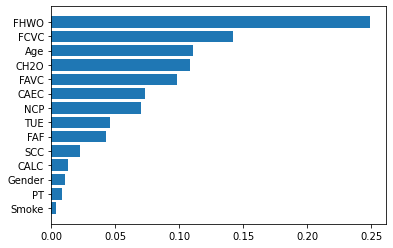}
    \caption{Random Forest Feature Importance Ranking.}
    \label{fig:2}
\end{figure}

The ranking of variable importance provided by this classifier is shown in Figure \ref{fig:2}. As one can observe, FHWO and Age are amongst the top three most important features in the model. This is a consistent factor as these are variables with high predictive power in classifying obesity. However, individuals cannot change their values on these variables. Therefore, in order to design strategies to combat obesity, they would have to design a ranking avoiding these variables. For instance, if we wanted to highlight the two most important variables in combating obesity that individuals can intervene in, based on this ranking, we would choose variables FCVC and CH2O. The problem with this ranking is that the fixed variables, although immutable, have an impact on the order of importance of the other variables, and this could lead to choosing suboptimal strategies for change. In the following sections, we illustrate how our models avoid this problem by creating rankings of variables according to the correct allocation of effort.

\subsection{Probability estimation}

In our approaches, the probabilities of taking the right or left branch, with or without effort, of all tree nodes are assumed to be known. However, in some situations, as in the case under analysis, these probabilities are unknown, and therefore must be estimated. In general, in datasets that follow a recorded evolution over time, one can make estimates based on these changes. However, when the data do not present this temporal character, this strategy is not possible and it is necessary to look for an alternative. For this purpose, we have designed a frequentist approach to estimate probabilities.

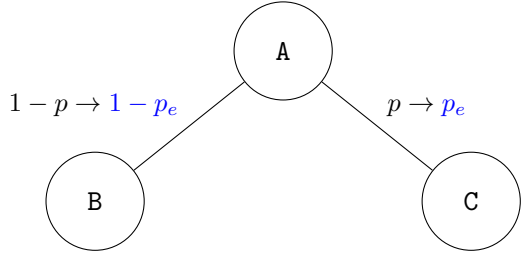
\begin{figure}[h]
    \centering

\begin{tikzpicture}[level distance=2cm,
level 1/.style={sibling distance=5cm}]

\node[circle,draw, inner sep=10] (Root) {\texttt{A}}
    child {
    node[circle,draw, inner sep=10](left1) {\texttt{B}} 
}
child {
    node[circle,draw, inner sep=10](right1) {\texttt{C}}
};
\path (Root) -- node[above right] {\small $p \rightarrow {\color{blue} p_e}$} (right1);
\path (Root) -- node[above left] {\small $1-p \rightarrow {\color{blue} 1-p_e}$} (left1);
\end{tikzpicture}
    \caption{Estimation node example}
    \label{fig:3}
\end{figure}

For each of the nodes, we will find a situation as in Figure \ref{fig:3}, where at node \texttt{A}, a split is applied to the a feature $X_j$ with the threshold value $c_j$ to reach either node \texttt{C} with probability $p$, or node \texttt{B} with probability $1-p$. In turn, these probabilities will be replaced by $p_e$ and $1-p_e$ when the individual applies effort in changing his variable $X_j$. 

To estimate the probabilities, we run $n_s = 1000$ simulations by perturbing the original value of the individual in the feature $X_j$, $x_j$, to a value $\hat{x_j}$, which will vary according to the type of variable $X_j$ and whether or not effort is applied:

\begin{itemize}
\item In case $X_j$ is a continuous or ordinal variable, we define $\hat{x_j}$ as $\hat{x_j} = x_j \pm \delta$, for a predefined value $\delta>0$ which is different in case some effort is applied or not:
\begin{itemize}
\item If no effort is applied, the positive or negative sign is determined by a Bernoulli distribution of probability $0.5$, and $\delta$ is calculated according to a uniform variable on the real segment from $0$ to $\sigma_j$, i.e. $\delta\sim U[0,\sigma_j]$, where $\sigma_j$ is the standard deviation of $X_j$ in the training sample. \item When some effort is applied, the sign is determined according to the favorable effect of $X_j$ with respect to the target variable (see Table \ref{tab:1}), and $\delta$ is allowed to take on larger values as $\delta\sim U[0,1.5\sigma_j]$.
\end{itemize}
\item In case $X_j$ is a binary feature:
\begin{itemize}
    \item If no effort is applied, we consider that an individual can switch between the two categories. We do this switching in the following way: we calculate the probability of belonging to the most represented class in the variable, $p_j$, and we change the value according to the realization of a $1-p_j$ probability Bernoulli variable. 
    \item If some effort is applied, we consider that the individual can only switch towards the beneficial class, and that it will not change if the value already is the beneficial one (for example, if a non-smoker person puts effort into no smoking, we assume that it will be achieved with probability $1$). As in the previous case, when dealing with the non-beneficial class, we calculate $p_j$ as the probability of belonging to the majority class, but in this case the switch is made according to a realization of a Bernoulli variable of parameter $\max(1-p, 0.2)$.
    \end{itemize}
\end{itemize}

After this, we count how many of the perturbed values take the right branch, $n_r$, and we define the probability as the ratio $n_r/n_s$.

\subsection{Rankings and results}

To this point, we have used a train set with 1838 individuals (2/3 of the data) to train a random forest of 25 trees of depth 5 with a $93\%$ accuracy on detecting obesity on the train set, and a $90\%$ of accuracy on the remaining 607 individuals of the test set. In addition, we have described the process of estimating probabilities for the individuals in the train  set on which we run our models. In this section, we describe how we apply our models, present some of the variable rankings generated and describe and report the results obtained in the test simulations.\\

In this experiment we consider that putting effort into more than $4$ variables is already impractical for an individual (i.e., parameter $\eta$ ranges in $\{1,2,3,4\}$). For each of the values of $\eta$, we run our three models on those individuals in the train set who belong to the obese category. Furthermore, for the $\kappa$-path problem, we explored three different versions varying the value of $\kappa$ from $33\%$, $50\%$ and $66\%$ of the leaf nodes, and taking $\mu = 10^{-6}$, which is a hyperparameter that can be calibrated. We set this value after some preliminary studies with good results, since being a small value, it leaves aside the percentage of very improbable solution paths.\\

For each of the individuals in the train set we ran all the different versions of our models, and for each solution obtained, we outlined the $\eta$ variables selected where the effort was placed. After this, we established an importance ranking based on the number of times the variable had been selected among the individuals in the sample. For the case of the random forest ranking (RFR), we used the results shown in Figure 2 and ignore the variables that individuals cannot modify. Thus, the ranking with $\eta=1$ would indicate using the variable FCVC, and with $\eta=4$ would indicate the variables FCVC, CH2O, FAVC and CAEC. Furthermore, we used as a control the average performance of three random selection rankings (RSR), where we simply randomly select $\eta$ variables from which individuals can change. \\

Once we had established the rankings and defined the variables on which the efforts are placed, we moved on to the evaluation phase. The evaluation was carried out on the 309 individuals of the test set who were obese. For each of the individuals, we performed 100 simulations with perturbed values in all variables following the same approach that was used to estimate the probabilities of change in the tree nodes, putting effort on the $\eta$ variables indicated by each ranking. After this, for each of the simulations, we re-used the classifier and calculate the percentage of individuals who were now classified as healthy. We then averaged the results of the $100$ simulations across all individuals and averaged the results to assess which selection of variables had a better performance. Table \ref{tab:2} shows the averaged results of the 100 simulations with the rankings of the models for the different values of $\eta$.\\

\begin{table}[h]
    \centering
    \begin{adjustbox}{width=\textwidth}
    \begin{tabular}{c|c|c|c|c|c|c|c}
    \toprule
         & RSR  & RFR & min-path & 
        $33\%$-path & $50\%$-path &
        $66\%$-path & max-path \\
        \midrule
         $\eta = 4$ & 25.22 & 26.07 & 26.85 & 24.62  & \textbf{31.81} & 24.91 & 24.01 \\
        $\eta = 3$ & 23.89 & 24.75 & 21.48 & 21.76 & \textbf{32.23} & 28.45 & 28.45 \\
        $\eta = 2$ & 25.10 & 18.20 & 20.91 &  \textbf{27.15} & 25.67 & 25.15 & 22.74 \\
        $\eta = 1$ & 21.80 & 18.32 & 21.97 & \textbf{26.66} & 26.53 & 21.89  & 21.67 \\
        \bottomrule
    \end{tabular}
    \end{adjustbox}
    \caption{Average  percentage of individuals reclassified as healthy based on our simulations for the different methodologies.}
    \label{tab:2}
\end{table}

Note that, among the obese individuals in the test sample, there are some who have very poor healthy habits, and for whom it would be practically impossible to change their classification (with any of the strategies). In order to show more realistic results, we normalize the information shown in the above table. To construct this table, we compute the percent values with respect to the total number of individuals which are feasible to change their class. The feasible-to-change sample size is obtained as the mean number of reclassified as healthy in 100 simulations by putting effort on all variables, also extending the value of the continuous variables to $1.5\sigma$ of the distribution in the favorable direction in all cases (instead of calculating a uniform value between 0 and $1.5\sigma$). Under this premise, we obtain a upper bound for the percent number of individuals that can change their class of 34.53$\%$ of the overall test sample. With this baseline, the results on the average number of individuals that adequately change their class from obese to healthy are shown in Table \ref{tab:4}.

\begin{table}[h]
    \centering
    \begin{adjustbox}{width=\textwidth}
    \begin{tabular}{c|c|c|c|c|c|c|c}
    \toprule
         & RSR  & RFR & min-path & 
        $33\%$-path & $50\%$-path &
        $66\%$-path & max-path \\
        \midrule
         $\eta = 4$ & 73.03 & 75.49 & 77.75 & 71.30 & \textbf{92.12} & 72.14 & 69.53 \\
        $\eta = 3$ & 69.18 & 71.67 & 62.20 & 63.01 & \textbf{93.33} & 82.39 & 82.39 \\
        $\eta = 2$ & 72.69 & 52.70 & 60.55 & \textbf{78.62} & 74.34 & 72.83 & 65.85 \\
        $\eta = 1$ & 63.13 & 53.05 & 63.62 & \textbf{77.20} & 76.83 & 63.39 & 62.75 \\ 
        \bottomrule
    \end{tabular}
    \end{adjustbox}
    \caption{Average  percentage of feasible-to-change individuals reclassified as healthy based on our simulations for the different methodologies.}
    \label{tab:4}
\end{table}

In both tables, we highlight in boldface those percent values that reach the maximum recovery for a given effort bound $\eta$.

Although the results on both tables are the same, in the second table one can observe that the ability to adequately re-classify the individuals is particularly impressive for the $\kappa$-path methodology, where in some cases, the percent of feasible-to-change individuals that were re-classified is more than $90\%$ by adequately selecting the features where putting effort.

Analyzing more in detail the results of the tables, in general we can see that the results when $\eta  = 3$ are better than those obtained when $\eta = 2$ and $\eta = 1$, and that they do not differ too much from those obtained when $\eta  = 4$. Another observation to note is that the best performing model is the $50\%$-path model, which proposes an intermediate optimization between maximizing the most favorable scenario and the least favorable scenario. On the other hand, when analyzing the solutions when $\eta = 1$, we observed that the FAVC variable had the greatest impact on reclassification individually. This variable was chosen in most of our models, and was also included in some of the random rankings. However, in RFR this variable is not included until we reach a value of $\eta = 3$, which partly explains the poorer performance of this ranking.

\begin{figure}[h]
    \centering
    \includegraphics[width=0.75\linewidth]{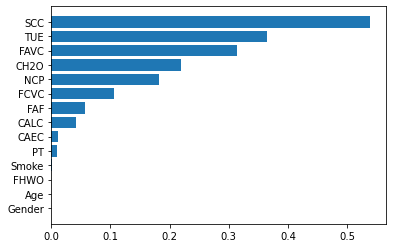}
    \caption{$\eta = 3, 50\%$-path Variable Importance }
    \label{fig:4}
\end{figure}

As a last detail, in Figure \ref{fig:4} we show the importance ranking of the variables obtained with the solution to the $50\%$-path problem  with $\eta=3$. According to this, we can see that the three top variables chosen for the ranking would be SCC (calories consumption monitoring), TUE (time using technology devices) and FAVC (frequent consumption of high caloric food). The case of the variable SCC is particularly curious, since only a few individuals in the sample monitor their daily calorie intake, and therefore it makes sense that this variable appears in the lower positions of the ranking provided by the Random Forest. This can be due to the fact that, as a binary variable with marginal classificatory impact on the total sample, it will be less likely to appear at higher levels of the ensemble trees, where variables with a higher impact on entropy will be preferred. 

The algorithmic artifacts relating to the importance of unbalanced binary variables that arise when using the Random Forest model may go unnoticed, or they may have a major impact as our computational results suggest (where for the case $\eta = 3$ our variable ranking works much better than the RFR in terms of success in guided future interventions---see Table \ref{tab:4}). In this sense, our methodology has a strong point in favour given that for an individual, all variables are a priori susceptible to be chosen for improvement.  Thus, our methods allow uncovering the features that, not being primarily relevant by their effect in classification based in the reduction of entropy of the total population, have a real capacity for reclassifying individuals.

\section{Conclusions}\label{sec:5}

In this paper, we provide a novel methodology to perturb the features of a given observation to be reclassified to a desired target class. The method integrates both the decisions on the features where the individual must put effort (among those affordable by a given budget), and the probabilities that the individual has to perturb its features. With all these in mind, we provide different mathematical optimization models to this end. The first one computes the features where putting effort to achieve the target class by maximizing the global probabilities. The second and third models are robust versions of the latter. In the second model the features are used to maximize the minimum probability among the possible paths to achieve the target class, whereas in the third model we apply a $k$-sum-based policy in which the perturbation with the probability of the $\kappa$ least probable paths is maximized. 

The practical validity of our approach is proved in a real dataset. As a result of each of these approaches, we get a ranking of the features providing information about the most viable features to achieve a class. In this study, we conclude that detecting the most (probabilistically) feasible features for an individual may have a big impact on the possibility of reaching the desired class. 

Possible extensions of our methodology include the study of these probabilistic approaches to other types of classifiers. The main difficulty of adapting the approach to other classification models is that it requires to keep track of the paths to reach a desired class, which is not always possible in black-box methods. It would be possible to apply this approach to multiclass classification methods, as the one in \cite{blanco2020optimal}. The study of non-independent features also deserves to be studied in the future, in particular its incorporation to the optimization models that we propose.

\section*{Acknowledgements}

This research has been partially supported by Spanish Ministerio de Ciencia e Innovación, AEI/FEDER grant number PID 2020 - 114594GBC21 and Junta de Andalucía, C-EXP-139-UGR23, and AT 21\_00032. The first author was also partially supported by the IMAG-Maria de Maeztu grant CEX2020-001105-M /AEI /10.13039/501100011033. This project is funded in part by Carnegie Mellon University’s Mobility21 National University Transportation Center, which is sponsored by the US Department of Transportation.


\end{document}